\input amstex
\documentstyle{amsppt}  \magnification=\magstep1  \overfullrule=0pt 
\NoRunningHeads \TagsOnRight \nologo  

\def\mb{\bold m} \def\ol{\overline} 
\def\Mb{\ol M} \def\mbgn{\ol M_{g,n}}  \def\mbn{\ol M_{0n}}
\def\w{\omega} \def\k{\w} \def\bw{\pmb\w} \def\wgn{\w_{g,n}}
\def\cn{\Cal C_n}  \def\wcm{\w_{\Cal C/M}}
\def\del{\partial}  \def\s{\sigma}
\def\rhs{r\.h\.s\. }  \def\lhs{l\.h\.s\. }
\def\Q{\Bbb Q}  \def\N{\Bbb N} \def\A{\Cal A} 
 \def\LR{\;\longleftrightarrow\;} \def\NIC{\,\text{\bf CohFT}_{\bold1}(k)\,}

\topmatter\title Higher Weil-Petersson Volumes of Moduli Spaces of Stable $n$-pointed Curves  \endtitle
  \author  R. Kaufmann, Yu. Manin, D. Zagier \endauthor
  \affil Max-Planck-Institut f\"ur Mathematik, Bonn, Germany \endaffil
  \abstract Moduli spaces of compact stable $n$-pointed curves carry a hierarchy 
of cohomology classes of top dimension which generalize the Weil-Petersson volume forms 
and constitute a version of Mumford classes. We give various new formulas for the integrals 
of these forms and their generating functions. \endabstract  \endtopmatter
\vskip -.5cm \line{\hfill{\it Dedicated to the memory of Claude Itzykson}}  \bigskip

\comment
\document
\null\bigskip
\centerline{\bf HIGHER WEIL--PETERSSON VOLUMES}
\centerline{\bf OF MODULI SPACES OF STABLE $n$--POINTED CURVES}
\medskip \centerline{\bf R. Kaufmann, Yu. Manin, D. Zagier}
\smallskip \centerline{Max--Planck--Institut f\"ur Mathematik, Bonn, Germany}
\bigskip  \line{\hfill{\it Dedicated to the memory of Claude Itzykson}}  \medskip

{\bf Abstract}. Moduli spaces of compact stable $n$-pointed curves carry a hierarchy 
of cohomology classes of top dimension which generalize the Weil-Petersson volume forms 
and constitute a version of Mumford classes. We give various new formulas for the integrals 
of these forms and their generating functions.
\endcomment

\bigskip  \centerline {\bf 0. Introduction} \medskip

Let $\mbgn$ be the moduli space of stable $n$-pointed curves of 
genus $g$. The intersection theory of these spaces is understood in the sense 
of orbifolds, or stacks. The algebro-geometric study of the Chow
ring of $\ol M_{g,0}$ was initiated by D. Mumford.

The following important version of Mumford classes on $\mbgn$ was introduced in [AC]. 
Let $p_n: \Cal {C}_n \rightarrow \mbgn$ be the universal curve, $x_i\subset\cn,i=1,\dots,n$,
the images of the structure sections, $\wcm$ the relative dualizing sheaf. Put for $a\geq 0$
  $$\w_n(a)=\w_{g,n}(a):=p_{n*}(c_1(\wcm(\sum_{i=1}^nx_i))^{a+1})\in H^{2a}(\mbgn,\Q)^{\Bbb S_n}\tag0.1$$
(We use here the notation of [KMK]; [AC] denote these classes $\kappa_i$. 
We will mostly omit $g$ in our notation but not $n$.)

The class $\w_{g,n}(1)$ is actually $\frac 1{2\pi^2} [v^{WP}_{g,n}]$ where 
$v^{WP}_{g,n}$ is the Weil-Petersson 2-form so that
   $$ \int_{\mbgn} \w_{g,n}(1)^{3g-3+n} =(2\pi^2)^{3g-3+n}\times \text{WP-volume of } \mbgn. \tag 0.2 $$
(see [AC], end of Sec.~1). Generally, we will call {\it higher WP-volumes} the integrals of the type
  $$ \int_{\mbgn} \wgn(1)^{m(1)} \dots \wgn(a)^{m(a)} \dots\,, \qquad \sum_{a\geq1}a m(a) = 3g-3+n. $$
The objective of this paper is to derive several formulas for these volumes and their generating functions.

For genus zero, we prove a recursive formula and a closed formula for them. The latter formula represents each
higher volume as an alternating sum of multinomial coefficients. It generalizes to the
higher genus, with the multinomial coefficients replaced by the correlation numbers 
$\langle \tau_{d_1} \dots \tau_{d_n} \rangle$ which are computable via Witten-Kontsevich's theorem [W], [K1]. 
In the genus zero  case, however, we can give an even better formula for the higher WP-volumes.
We first encode these values via a generating function in infinitely many variables and translate the
recursions for them into an infinite system of non-linear differential equations for this generating
function. It then turns out that the inverse power series of (a slightly modified version of) this generating 
function satisfies a system of {\it linear} differential equations, which can then be solved explicitly.

We will now explain our results for the case of the classical WP-volumes (0.2) of the genus
zero moduli spaces first calculated by P.~Zograf [Z]. Put $v_3 = 1$ and
  $$ v_n := \int_{\mbn} \w_n (1)^{n-3}, \qquad n \geq 4\,. \tag 0.3$$
The main result of [Z] is:
  $$ v_n= \frac 12\sum_{i=1}^{n-3}\frac{i(n-i-2)}{n-1}
  {{n-4}\choose{i-1}}{{n}\choose{i+1}}\,v_{i+2}\,v_{n-i}\,,\qquad n\geq 4\,.\tag 0.4 $$
(The factor 1/2 was inadvertently omitted when this formula was quoted in [KMK].)  Consider the generating functions
  $$ \Phi(x)=\sum_{n=3}^\infty\frac{v_n}{n!(n-3)!}\,x^n\,,\quad 
        h(x)=\Phi'(x)=\sum_{n=3}^\infty\frac{v_n}{(n-1)!(n-3)!}\,x^{n-1}\,.\tag 0.5 $$
Then (0.4) directly translates into the differential equation
  $$ xh''-h' = (xh'-h)\,h''\,.\tag0.6$$
Notice that according to [Ma], the generating series (0.5) arise in the Liouville gravity models.

In the first part of our paper we generalize both (0.4) and (0.5) to arbitrary WP-volumes of genus zero: see
Theorems 1.2.1 and 1.6.1. We follow the scheme of proof explained in [KMK], sec. 3.  The second series of our  
results gives more explicit formulas which specialize to the case of $v_n$ in the following way: for $n \geq4$,
  $$ v_n=\sum_{k=1}^{n-3}\frac{(-1)^{n-3-k}}{k!}\sum\Sb m_1,\ldots,m_k>0\\m_1+\dots+m_k=n-3\endSb
   \binom{n-3}{m_1,\dots,m_k}\binom{n-3+k}{m_1+1,\dots,m_k+1}\,.\tag 0.7 $$
(Actually, we prove an analog of (0.7) for arbitrary genus, but as we have already mentioned the multinomial 
coefficients are then replaced by the less well understood numbers $\langle \tau_{d_1} \cdots \tau_{d_n} \rangle$.)
Then, using either this explicit formula or by inverting the system of differential equations, we obtain
a formula for the generating function of the WP-volumes as the inverse power series of a very simple
power series.  For the case of the $v_n$ this goes as follows:
Differentiating the differential equation (0.6) gives $h'h'''=x{h''}^3$ or (setting $y=h'$) $yy''=x{y'}^3$. 
This is now cubic rather than quadratic, but if we interchange the roles of $x$ and $y$ then it miraculously 
transforms into the Bessel equation $y\,\dfrac{d^2x}{dy^2}+x=0$, leading to the explicit solution of (0.4)
via an inverted modified Bessel function:
  $$ y = \sum_{n=3}^\infty\frac{v_n}{(n-2)!\,(n-3)!}\,x^{n-2} \quad\Longleftrightarrow\quad
    x = \sum_{m=1}^\infty \frac {(-1)^{m-1}}{m!\,(m-1)!}\,y^m\;. \tag 0.8 $$
It is tempting to see this as another tiny bit of the general ``mirror phenomenon" first observed
for Calabi-Yau threefolds.

As a corollary of (0.8) we get the asymptotics
  $$ v_{n+3} \approx \frac{(2n)!}{{\roman C}^n} \left (1.3620537\ldots 
  - \frac{0.131\dots}{n} + \frac {0.019\dots}{n^2} - \dots \right ) \tag 0.9 $$
where $\roman C = 2.496918339\dots$ is a constant that can be expressed in terms of Bessel functions and their
derivatives.  The existence of such an asymptotic formula---for all genera, and with a constant $C$ independent
of the genus---was conjectured by Claude Itzykson (P.~Zograf, private communication).

It is interesting to compare (0.9) to the asymptotics of the Euler characteristic
  $$\chi(\ol M_{0,n+3})\approx\frac1{\sqrt{n+2}}\left(\frac{n+2}{e^2-2e}\right)^{n+\frac32} \tag0.10 $$
(cf. [M], p. 403). One more problem in the same spirit is to study the asymptotic structure of the representation 
of $\Bbb S_n$ on $H^*(\mbn)$ with respect to the Plancherel measure, in the same sense as it was done for the 
regular representation in [LS] (cf. also [VK].) A relevant generating function was recently calculated by 
E.~Getzler and M.~Kapranov, cf. [G].

The paper is structured as follows. In \S1 we prove the recursive relations 
and the differential equations for the generating function in genus zero. 
In \S2 we derive by a different method explicit formulas for higher WP-volumes and prove
the analogue of the inversion formula (0.8).
Finally, \S3 briefly explains our main motivation for studying WP-volumes: they and their generating 
function naturally arise in the theory of the so called Associativity Equations and the operation of tensor 
product in the category of algebras over the $\{H_*(\ol M_{0,n+1})\}$-operad. (For more details see [KMK].)
Finally, in a short appendix we make some remarks about the asymptotic formulas (0.9) and (0.10)
 and about the Betti numbers of the spaces $\ol M_{0,n+3}$.

\newpage
\centerline{\bf \S \, 1. Recursive relations and }
\centerline{\bf differential equations for the generating function}
\bigskip

{\bf 1.1. Notation.} Consider the semigroup $N^\infty$ of sequences $\mb=(m(1),m(2),\dots)$ where
$m(a)$ are nonnegative integers, $m(a)= 0$ for sufficiently large $a$.

We put 
  $$V_g (\mb) := \frac1{(\sum_{a\geq 1} a m(a))!} \int_{\mbgn} \prod_{a \geq 1}
    \frac {\w_{g,n}(a)^{m(a)}}{m(a)!} \in \Q \tag 1.1 $$
interpreting the \rhs as zero unless $\sum_{a\geq 1}am(a)=\roman{dim}\,\mbgn=3g-3+n$.
In the rest of this section $g=0$ and we write $V(\mb)$ instead of $V_0 (\mb)$. 

In shorter versions of expressions like (1.1) we will use notation of the type
  $$\aligned &|\mb|:=\sum_{a\geq1}am(a),\quad\|\mb\|:=\sum_{a\geq1}m(a),\quad\mb!:=\prod_{a\geq1}m(a)!,\\
   &\qquad\bw_n^\mb=\prod_{a\geq1}\wgn(a)^{m(a)},\quad \bold s^\mb=\prod_{a\geq1}s_a^{m(a)}\endaligned\tag1.2$$
where $\bold s=(s_1,s_2,\dots)$ is a family of independent formal variables or complex numbers. For
instance, we have $V(\mb)=\int\bw^\mb/\mb!\,|\mb|!\,$ in this notation.

\smallskip
{\bf 1.2. A recursive Formula for $V(\mb)$.} Put
  $$K(n_1,\dots,n_a):=\frac1{n_1(n_1+n_2)\cdots(n_1+\dots+n_a)} \tag 1.3$$
and denote by $\delta_a \in N^\infty$ the sequence with 1 at the $a$-th place and zeros elsewhere.

\proclaim{\quad 1.2.1. Theorem}  For any $\mb$ and $a\geq1$, we have:
  $$ (m(a)+1) V(\mb+\delta_a) = (|\mb|+a+1) \sum_{\mb = \sum_{i=1}^{a+1}\mb_i} 
   K(n_1, \dots ,n_a) \prod_{i=1}^{a+1} V(\mb_i) \tag 1.4 $$
where in each summand of (1.4)
  $$(n_1, \dots ,n_a) := (|\mb_1|, \dots |\mb_a|) + (2,1,\dots ,1) \tag 1.5 $$
{\rm (notice the absence of $|\mb_{a+1}|$)}. These relations uniquely define $V(\mb)$
starting with $V({0}) = 1$.  \endproclaim

{\bf 1.2.2. A particular case of (1.4).}
Applying (1.4) to $V(m):= V(m,0,0,\dots)$ and $a=1$ we get:
  $$\split (m+1)V(m+1) &= (m+2) \sum_{m=m_1+m_2} \frac 1{m_1+2} V(m_1) V(m_2)\\
   &= \frac12(m+2)\sum_{m=m_1+m_2}(\frac 1{m_1+2}+\frac 1{m_2+2})\,V(m_1)\,V(m_2)\endsplit $$
so that
  $$ \frac{V(m+1)}{m+3} = \frac12\,\frac{(m+2)(m+4)}{(m+1)(m+3)}\,
   \sum_{m=m_1+m_2} \frac {V(m_1)}{m_1+2}\,\frac {V(m_2)}{m_2+2}.$$
On the other hand, Zograf's recursive relations (0.4) can be rewritten as 
  $$ \frac{(n-2)v_n}{(n-3)!(n-1)!}=\frac12\frac{(n-2)n}{(n-3)(n-1)}\sum\Sb n+2=p+q\\p,q\ge3\endSb 
     \frac{(p-2)v_p}{(p-3)!(p-1)!}\,\frac{(q-2)v_q}{(q-3)!(q-1)!}. $$
These relations agree for $V(n-3)={v_n}/{(n-3)!^2}$ which is the correct formula in view of (0.3) and (1.1).

\smallskip
{\bf 1.2.3. Another special case.} Let us use (1.4) to compute $V(\mb)$ for $|\mb|\le2$. 
We have from (1.4) for $\mb =0$:
  $$ V(\delta_a) = (a+1)K(2,1,\dots,1) = \frac1{a!}. \tag 1.6 $$
But already in the next case, the consistency of the two possible formulas obtained from (1.4) 
for $V(\delta_a +\delta_b)$ is not evident a priori.

Put $\mb=\delta_b$. There are $a+1$ partitions of $\delta_b$ contributing to (1.4):
in the $k$-th partition $\mb_k=\delta_b,\;\mb_j=0$ for $j\ne k,\;\prod_{b=1}^{a+1} V(\mb_i)=1/{b!}$,
  $$ (n_1,\dots,n_a)=\cases(2,1,\dots,1)+b\delta_k &\text{for $k\leq a,$}\\ (2,1,\dots,1) &\text{for $k=a+1$,}\endcases$$
so that
  $$ K(n_1, \dots ,n_a) = \frac1{(a+b+1)!} \frac{(b+k)!}{k!} $$
and
  $$ 2V(\delta_a+\delta_b)=(a+b+1)\sum_{k=1}^{a+1}\frac{(b+k)!}{(a+b+1)!k!}
    \frac1{b!}=\frac1{(a+b)!}\sum_{k=1}^{a+1}\binom{b+k}k$$
This does not look symmetric in $a,b$, but of course it is:
  $$ V(\delta_a+\delta_b) = \frac 1{2(a+b)!}\left[\binom{a+b+2}{a+1}-1\right]\tag 1.7$$
We will generalize (1.6) and (1.7) below to all $n$ (and $g$).

\bigskip
We will now start proving Theorem 1.2.1. We will use some of the formalism of 
[KMK]. For any stable $n$-tree $\s$, denote by 
$\varphi_\s:\ol M_\s \hookrightarrow \mbn$  the embedding of the 
corresponding closed stratum. We recall that the image of a generic point of 
$\ol M_\s$ parametrizes a curve whose dual tree is $\s$. The set of 
vertices $v \in V_\s$ of this tree bijectively corresponds to the set of 
irreducible components of the curve. The edges $e \in E_\s$ ``are'' 
singular points of the curve, and the unpaired flags (tails or leaves) are in a 
bijection with the labelled points, that is, with $\{1, \dots ,n \}$. 
For $v \in V_\s$, we denote by $|v|$ the number of flags incident to $v$.

For any stable $n$-tree $\s$, we put (with notation (1.2))
  $$ \Omega_n(\mb,\s) = \int_{\ol M_\s} \frac {\varphi^*_\s(\bw_n^\mb)}{\mb !}, \tag 1.8 $$
interpreting this as zero unless 
$n-3-|\mb| = \roman{codim}\, \varphi_\s(\ol M_\s) = |E_\s|$.
If $\s_n$ is an one-vertex $n$-tree, we write 
$\Omega_n(\mb):=\Omega_n(\mb,\s_n).$
Notice that $\Omega_n(a)$ from [KMK] is 
$\left(\frac{n-3}{a}\right)!\,\Omega_n(\frac{n-3}{a}\delta_a)$
in our present notation. The numbers $V(\mb)$ in (1.1) are 
$\Omega_n(\mb)/|\mb|!$.

\proclaim{\quad 1.3. Lemma} We have 
  $$ \Omega_n(\mb,\s) = \sum\Sb(\mb_v|v\in V_\s):\\ \mb=\sum\mb_v\endSb\,
  \prod_{v\in V_\s}\Omega_{|v|}(\mb_v),\tag1.9 $$
where the sum in \rhs is taken over all partitions of $\mb$ indexed by vertices of $\s$. \endproclaim

{\bf Proof.} This follows from the crucial fact that $\w_n (a)$
form what is called a ``logarithmic CohFT'' in [KMK], sec. 3,
i\.e\. satisfy the additivity property established in [AC], (1.8), for any genus:
  $$ \varphi_\s^*(\w_n(a)) = \sum_{w\in V_\s} \roman{pr}^*_w(\w_{|w|}(a)) \tag 1.10$$
where we identify $\Mb_\s$ with $\prod_{w\in V_\s} \Mb_{0,|w|}$ and $\roman{pr}_w$ 
means the respective projection. (Notice that although these identifications are defined only 
up to the action of $\prod_{w\in V_\s} \Bbb S_{|w|}$, the classes $\roman{pr}^*_w(\w_{|w|}(a))$
do not depend on their choice, being $\Bbb S_{|w|}$-invariant).

Hence
$$
\split
\int_{\Mb_\s} \varphi_\s^*(\bw_n^\mb) 
 &= \int_{\prod_{v\in V_\s} \Mb_{0,|v|}} 
   \prod_{a \geq 1} 
     \left( \sum_{w\in V_\s} \roman{pr}^*_w(\w_{|w|}(a))\right )^{m(a)}\\
 &=  \int_{\prod_{v\in V_\s} \Mb_{0,|v|}} 
   \prod_{a \geq 1} 
   \sum \Sb (m_w(a)|w \in V_\s): \\ m(a) = \sum_w m_w(a) \endSb
  \frac {m(a)!}{\prod_w m_w(a)!} 
  \prod_{w\in V_\s} \roman{pr}^*_w(\w_{|w|}(a))^{m_w(a)}\\
 &= \int_{\prod_{v\in V_\s} \Mb_{0,|v|}} 
   \sum \Sb (\mb_w|w \in V_\s): \\ \mb = \sum_w \mb_w \endSb
  \frac {\mb!}{\prod_w \mb_w!}
  \prod_{w\in V_\s} \roman{pr}^*_w(\bw_{|w|})^{\mb_w} \\
 &= \sum \Sb (\mb_w|w \in V_\s): \\ \mb = \sum_w \mb_w \endSb
 \mb! \prod_{v\in V_\s} \int_{ \Mb_{0,|v|}} \frac {\bw_{|v|}^{\mb_v}}{\mb_v!}
\endsplit
 \tag Q.E.D.
$$

{\bf 1.4. Calculation of $\w_n(a)$ via strata classes.}
For a fixed $n \geq 3$ and $a \geq 1$ consider labelled $(a+1)$-partitions
  $$ S:\quad \underline n := \{1,\dots,n\} = S_1 \amalg \dots \amalg S_{a+1}. $$
Denote by $\tau(S)$ the $n$-tree with $V_{\tau(S)} = \{v_1, \dots , v_{a+1} \}$,
and edges connecting $v_i$ to $v_{i+1}$ for $i=1,\dots,a$, and unpaired flags
(numbered by) $S_i$ put at $v_i$. The stability condition for
$\tau(S)$ and $S$ is:
  $$ n_i:=|S_i|\geq2\;\text{for $i=1,a+1$}; \quad\geq1\;\text{for $i=2,\dots,a$}. \tag 1.11 $$
In the following proof, all partitions are stable. Denote by $m(S)$ the dual 
cohomology class of the cycle $\varphi_{\tau(S)}(\Mb_{\tau(S)})$ in $\mbn$.

\proclaim {\quad 1.4.1. Lemma} We have
$$ \w_n(a) = \sum_{S: \underline{n} = S_1 \amalg \dots \amalg S_{a+1}}
  \frac {(n_1-1)(n_{a+1}-1)n_1 \dots n_{a+1}}{n(n-1)} K(n_1,\dots,n_a)\,m(S) \tag 1.12 $$
where $K(n_1,\dots,n_a)$ is defined in $(1.3)$. \endproclaim

{\bf 1.4.2. Notation.} To state some intermediate formulas we will need some of the notation of [KMK]. 
Let $T_n(a)$ be the set of $n$-trees with $a$ edges.
For any flag $f$ denote by $\beta (f)$ the set of tails of the branch of $f$ and $S(f)$ the set of their labels. 
Then to any set of flags $T$ we associate the set $S(T) := \bigcup_{f \in T} S(f) \subset \{1, \ldots , n\}$. 
If $\{S(T_1),S(T_2)\}$ is a partition of $\{1,\ldots,n\}$ we use the shorthand notation $D_{T_1,T_2}$ for 
$D_{S(T_1),S(T_2)}$.  Let $\tau$ be an $n$-tree and let $e$ be one of its edges, 
$\del e = \{v_1,v_2\}, \s_e$ the corresponding partition 
$S_1 \amalg S_2$ and $D_e$ the corresponding divisor. 
Choosing flags $\{i,j\} \in F(v_1)$ and $\{j,k\} \in F(v_2)$ we have [KM] the following formula:
  $$ \split    D_e m(\tau) = 
&-\sum_{{T:\,\{i,j\}\subset T\subset F(v_1)\atop 2\le |T|\le |F(v_1)|-1}}D_{T,F(v_2)\amalg F(v_1)\setminus T}m(\tau) \\
&-\sum_{{T:\,\{k,l\}\subset T\subset F(v_2)\atop 2\le |T|\le |F(v_2)|-1}}D_{T,F(v_1)\amalg F(v_2)\setminus T}m(\tau) 
 \endsplit \tag 1.13 $$

\proclaim{\quad 1.4.3. Definition} A tree is called linear if each vertex has at most two incident edges. 
An orientation of a linear tree is a labelling of its vertices by $\{1,\dots,|V(\tau)|\}$ such that $v_i$ 
and $v_{i+1}$ are connected by an edge for $i= 1, \dots, |E(\tau)|$. \endproclaim  

We denote by $LT_n(a)$ the set of stable linear $n$-trees with $a$ edges modulo isomorhism. Given
a geometrically oriented linear tree we number its vertices in the positive direction. 
\smallskip

{\bf 1.4.4. Remark.}
The oriented linear trees in $LT_n(a)$ are in 1-1 correspondence with labeled $a+1$-partitions
$S: \underline{n} := \{1,\dots,n\} = S_1 \amalg \dots \amalg S_{a+1}$ which satisfy (1.11).
\smallskip

{\bf 1.4.5. Tautological classes and the $\w_n(a)$.}
In the proof of the Lemma we will consider some additional classes in 
$H^*(\mbgn, \Q)$. Let $\xi_i: \mbgn \rightarrow \Cal C_n$
be the structure sections of the universal curve. Put as in [AC]:
  $$ \Psi_{n,i} := \xi_i^* (c_1(\wcm)) \in H^2(\mbgn, \Q). \tag 1.14 $$
Here we will need them only for $g=0$; see \S2 for any genus.

Identify $C \rightarrow \mbn$ with the forgetful morphism  $p_{n}: {\ol M}_{0,n+1} \rightarrow \mbn$. 
Then $\xi_i(\mbn)$ becomes the divisor $D_i=D_{\{i,n+1\}\{1,\ldots,\hat i,\ldots,n\}}$ in $\ol M_{0,n+1}$ and
  $$ \Psi_{n+1,i} = \varphi^*_{D_i}(-D_i^2) $$
where $\varphi^*_{D_i}$ denotes the pullback onto the divisor $D_i$.  We know from [AC]:
  $$ \w_{n-1}(a) = p_{n-1*} (\Psi^{a+1}_{n,n}) $$
Combining these two formulas we obtain:
  $$ \w_{n-1}(a) = p_{n-1*}\circ \varphi^*((-1)^{a+1} D_i^{a+2}) \tag 1.15 $$

To derive (1.15) notice that 
$$ \Psi_{n,i}=\sum_{n\in S\subset\{1,\ldots,n\}}\frac{(n-|S|)(n-|S|-1)}{(n-1)(n-2)}D_{S,\{1,\ldots,n\}\setminus S} $$
(see [KMK]) so that we have
$$ \split  \Psi^{a+1}_{n,n}
  &= \left (\sum_{n\in S\in\{1,\ldots,n\}}\frac{(n-|S|)(n-|S|-1)}{(n-1)(n-2)}D_{S,\{1,\ldots,n\}\setminus S}\right)^{a+1} \\
  &=\varphi^*_{D_n}(\left(\sum_{\{n,n+1\}\subset S\in\{1,\ldots,n+1\}}\frac{(n-|S|)(n-|S|-1)}{(n-1)(n-2)}     
     D_{S,\{1,\ldots,n+1\}\setminus S} \right )^{a+1} D_n)  \\
  &= \varphi^*_{D_n}((-1)^{a+1} D_n^{a+2}) \endsplit $$

{\bf 1.4.6. A calculation.} Denote by $D_nLT(k)$ the set of oriented linear $(n+1)$-trees with $a$ edges, 
whose monomials are divisible by $D_n$ and whose orientation is given by calling $v_{k+1}$ the trivalent vertex 
with the two tails $n$ and $n+1$.  Also take $S$ to be the set of the flags of the other vertex $v_{k}$ of the edge corresponding to $D_n$ without the flag belonging to that edge. Then
$$ D_n^k=\sum_{\tau \in LDT_n(a)}  \frac {|v_1|(|v_1|-1)}{(n-1)(n-2)} \prod_{i=2}^{k-1} 
\frac {|v_i|-2}{\sum_{j=i}^{k} (|v_j|-2)} m(\tau). \tag 1.16 $$

We will prove (1.16) by induction using the following versions of (1.13).
Let $\tau$ be a tree, which has an edge  $e$ corresponding to $D_n$, then call $v_2$ the vertex with 
$F(v_2)= \{n,n+1,f_e\}$, where $f_e$ is the flag corresponding to $e$.

Averaging the formula (1.13) over the set $S$ of all flags of $v_1$ without the flag belonging to $e$ we obtain:
  $$ D_n m(\tau) = - \sum_{{T\subset S \atop 2\leq |T|\leq |F(v1)|-2}}
  \frac{|T|(|T|-1)}{|S|(|S|-1)} D_{T, \{n,n+1\} \amalg S\setminus T} m(\tau)  \tag 1.17 $$

Fixing {\it one} particular flag $f$ of $S$ and averaging over the rest we obtain:
  $$ D_n m(\tau) = - \sum_{{f \in T\subset S \atop 2\leq |T|\leq |F(v1)|-2}}
  \frac {|T|-1}{|S|-1} D_{T, \{n,n+1\} \amalg S\setminus T} m(\tau) \tag 1.18 $$

Now for $k=1$ the formula (1.16) is clear and for $k=2$ it is a consequence of (1.18). For $k>2$ We have
  $$\align D_n^k 
   & = D_n\,D_n^{k-1} = D_n\,\sum_{\tau\in D_nLT(k-1)} \frac{|v_1|(|v_1|-1)}{(n-1)(n-2)}
    \prod_{i=2}^{k-2} \frac{|v_i|-2}{\sum_{j=i}^{k-1} (|v_j|-2)} m(\tau) \\
   & =\sum_{\tau \in D_nLT(k-1)} \sum_{f \in T \subset S} \frac{|v_1|(|v_1|-1)}{(n-1)(n-2)}\times\\ 
   &\qquad\quad \prod_{i=2}^{k-2}\frac {|v_i|-2}{\sum_{j=i}^{k-1} (|v_j|-2)} 
     \frac{|T|-1}{|v_{k-1}|-2} D_{T,\{n,n+1\}\amalg S\setminus T} m(\tau) \tag 1.19 \endalign $$
where we have used (1.18) with the distinguished flag being the unique flag of $S$ belonging to an edge. 
This guarantees that the sum will again run over linear trees.  In the second sum there is one edge inserted 
at the vertex $v_{k-1}$ giving two new vertices $v', v''$ with $|v'|+|v''|=|v_{k-1}|+2$.  
Giving $v', v''$ the labels $k-1$ and $k$ 
and labelling the old vertex $v_k$ with $k+1$ in the second sum we obtain the desired result (1.16).
\smallskip

{\bf 1.4.7. Proof of Lemma.}
What remains to be calculated is $p_{n-1*}\circ \varphi^*_{D_n}$ of the above 
formula for $D_n^{a+2}$. The only nonzero contributions come from trees 
$\tau \in D_nLT(a+2)$ with $|v_k|=3$, so that exactly one of the flags is
a tail. Hence after push forward and pullback the sum will run over 
oriented linear trees with the induced orientation given by the image of 
$v_i$ with a distinguished flag at the vertex $v_{k-2}$. 
Summing first over the possible distinguished flags amounts to multiplication by $|v_{k-1}|$.
We obtain:
$$
\w_n(a) = \sum_{\roman {oriented}\, \tau \in LT_n(a)} 
\frac {(|v_{a+1}|-1)(|v_1|-1)}{n} \prod_{i=1}^{a+1} 
\frac {|v_i|-2}{\sum_{j=i}^{a+1} (|v_j|-2) + 1} m(\tau)
$$
which using remark 1.4.4 can be rewritten as a sum over partitions 
$$
\w_n(a) =\sum_{S: S_1\amalg \dots \amalg S_{a+1}} \frac{n_1 n_{a+1}}{n}
(n_1-1)n_2 \cdots n_a (n_{a+1}-1) \frac 1{n-1} K(n_a+1,\ldots,n_2) m(S)
$$
with $n_i=|v_i|-1$ for $i=1,a+1$ and $n_i=|v_i|-2$ for $i=2,\dots,a$, which is equivalent to (1.12).
\smallskip

{\bf 1.4.8. Remark.} Instead of using (1.18) in the induction one can successively apply (1.17).
In this case one obtains a formula for $\w_n(a)$ involving all boundary
strata. Since not necessarily linear trees cannot be handled using only
partitions the associated generating functions and recursion relations become very complicated.
\smallskip

{\bf 1.5. Proof of Theorem 1.2.1.}  In view of (1.8), we have
  $$  \Omega_n(\mb+\delta_a) = \int_{\mbn} \prod_{b\geq 1} 
  \frac {\w_n(b)^{m(b)}}{m(b)!} \wedge \frac {\w_n(a)}{m(a)+1} \tag 1.20 $$
Instead of wedge multiplying by $\w_n(a)$ we can integrate the product ${\w_n^\mb}/{\mb!}$ over the cycle 
obtained by replacing $m(S)$ by $\varphi_{\tau(S)}(\Mb_{\tau(S)})$ in the \rhs of (1.12).
The separate summands then can be calculated using (1.8) and (1.9).  The net result is:
 $$  \multline  (m(a)+1)\Omega_n(\mb+\delta_a) = 
   \sum_{S: \underline{n}  = S_1 \amalg \dots \amalg S_{a+1}}
    \frac {(n_1-1)(n_{a+1}-1)n_1 \dots n_{a+1}}{n(n-1)}\times\\
 K(n_1,\dots,n_a) \sum_{\mb = \mb_1+ \dots+\mb_{a+1}}\Omega_{n_1+1}(\mb_1) \Omega_{n_{a+1}+1} (\mb_{a+1})
  \prod_{i=2}^a \Omega_{n_i+2}(\mb_i) \endmultline \tag 1.21 $$
Now, the product of $\Omega$'s vanishes unless
  $$|\mb_i|=n_i-2\;\text{for $i=1,\,a+1$}, \quad|\mb_i|=n_i-1\quad\text{for $i=2,\dots,a$} \tag 1.22 $$
so that $n=|\mb+\delta_a|+3$. Hence we can make the exterior summation over vector $(a+1)$-partitions of $\mb$, 
and for a fixed $(\mb_i)$ sum over the set of $(a+1)$-partitions  of $\underline{n}$ satisfying (1.22).
Since the coefficients in (1.20) depend only on $(n_i)$ rather than $(S_i)$, we can then replace the summation 
over $(S_i)$'s by multiplication by $\frac{n!}{n_1!\dots n_{a+1}!}\,$.  This leads to 
  $$ (m(a)+1)\frac{\Omega_n(\mb+\delta_a)}{|\mb+\delta_a|!} =(n-2)\!\!\sum_{\mb=\mb_1+\dots+\mb_{a+1}}
   \! K(n_1,\dots,n_a)\,\prod_{i=1}^{a+1}\frac {\Omega_{|\mb_i|+3}(\mb_i)}{|\mb_i|!} $$
which is equivalent to (1.4) in view of (1.8) and (1.1). Q.E.D.
\smallskip

{\bf Remark.} We do not know whether for $g\geq 1$ the classes $\w_n(a)$
belong to the span of the boundary strata and if yes, what might be a 
generalization of (1.20). Therefore we cannot extend the recursive relations (1.4) to arbitrary genus.
\smallskip

{\bf 1.6. The differential equation for a generating function.}  Put 
  $$ F(x;\bold s) = F(x;s_1,s_2,\dots ) := \sum_\mb V(\mb)x^{|\mb|}\,\bold s^\mb \in\Q[\bold s][[x]] \tag 1.23 $$
and denote $\del_a=\del /\del s_a$, $\del_x=\del/\del x$. Then the 
recursion (1.4) is equivalent to:

\proclaim {\quad 1.6.1. Theorem} $F$ satisfies the following system of differential equations:
  $$ \del_a F = \del_x\left(\sum_{k=0}^a(-1)^k\frac{F^{2k+1}}{({\del_x F})^{k+1}}
  \del_{a-k}F\right),\qquad a=1,2,\dots \tag 1.24 $$
where we put $\del_0=x\,\del_x$. It is the unique solution of this system in $1+x\Q[\bold s][[x]]$
with $F(x;\bold 0)=0$. \endproclaim

{\bf Proof.} Put $H_0(x;\bold s) = x$ and
  $$ H_a(x;\bold s):=x\sum_{\mb_1,\dots,\mb_a}K(n_1,\dots,n_a)\prod_{i=1}^a V(\mb_i)x^{|\mb_i|+1}\,\bold s^{\mb_i} $$
 for $a=1,\,2,\ldots$, where the summation is over $a$-tuples of vectors $\mb_i\in N^\infty$.
In particular, we have $H_1(x;\bold s) = \int_0^x \xi F(\xi;\bold s) d\xi$.
Multiply (1.4) by $ x^{|\mb + \delta_a|}\bold s^\mb$ and sum over all $\mb$. Taking into account that in each summand 
  $$ \split n_1+ \dots +n_a = |\mb|+a+1 &= |\mb + \delta_a| +1, \\
  (n_1+ \dots +n_a)K(n_1, \dots ,n_a) &= K(n_1, \dots ,n_{a-1}), \endsplit $$
we obtain
  $$ \del_a F(x;\bold s)= \del_x\bigl(H_a(x;\bold s)\,F(x;\bold s)\bigr), \qquad a\geq 1. \tag 1.25 $$
A similar calculation shows that
  $$ \del_x H_a(x;\bold s) = H_{a-1}(x;\bold s)\,F(x;\bold s), \qquad a \geq 1. \tag 1.26 $$
Combining these two identities, we get
  $$ \split \del_a F &= \del_x H_a \cdot F + H_a \cdot \del_x F\\
   &= H_{a-1} \cdot F^2 + H_a \cdot \del_x F, \qquad a\ge 1 \endsplit $$
or
  $$ H_a = \frac {\del_a F}{\del_x F} - H_{a-1}\frac {F^2}{\del_xF}\,. \tag1.27 $$
By induction starting with $H_0=x$ we obtain from here
  $$ H_a = \sum_{k=0}^a (-1)^k \,\frac{F^{2k}\,\del_{a-k}F}{(\del_xF)^{k+1}} \tag 1.28 $$
(recall that $ \del_0 := x \del_x$).  Substituting this into (1.25) gives (1.24). 

For the uniqueness, we reverse the argument.  Suppose that $F(x)=F(x;\bold s)$ satisfies (1.24) and
define $H_a$ by (1.28).  Then $H_0(x)=x$ (by the definition of $\del_0F$) and $H_a$ satisfies
(1.27) for $a\ge1$, while (1.24) says that $\del_aF=\del_x(H_aF)$.  Combining these equations gives
(1.25).  By assumption $F(x;\bold s)$ has a Taylor series $\sum_{n=0}^\infty A_n(\bold s)x^n$ where
$A_0(\bold s)\equiv1$ and $A_n(\bold s)$ for $n\ge1$ is a polynomial with no constant term. 
Equation (1.28) then shows that $H_a(0)=0$,
after which the equation $\del_x(H_a)=H_{a-1}F$ gives inductively $H_a=\frac{x^{a+1}}{(a+1)!}+\cdots$
where the coefficient of $x^{a+n+1}$ is a weighted homogeneous polynomial in $A_1,\ldots,A_n$ of
weight $n$.  The equation  $\del_aF=\del_x(H_aF)$ then gives a formula
for all derivatives $\del_aA_n(\bold s)$ ($a=1,2,\ldots$) as polynomials in $A_1(\bold s),\ldots,A_{n-1}(\bold s)$.
This fixes $A_n(\bold s)$ inductively up to a constant which is uniquely determined by the normalizing
condition $A_n(\bold 0)=0$. In this argument we have implicitly assumed that $\del_xF$ is an invertible power
series (i.e. that $A_1(\bold s)$ is not identically 0) in order to make sense of (1.24) in the ring of power
series, but what we have really proved does not need this assumption, namely, that $(F,H_a)$ is the 
unique solution in elements of $\Q[\bold s][[x]]$ of the system of differential equations (1.25), (1.26) 
subject to the normalizing conditions $F(0,\bold s)=F(x;\bold 0)=1$, $H_0(x)=x$, $H_a(0)=0$.
  
\smallskip
{\bf 1.6.2. Example.} In the special case $\bold s=(s,0,0,\ldots)$, the function (1.23) reduces to $F(x;\bold s)=f(xs)$
with $f(x)=\sum_{m\ge0}V(m)x^m$, $V(m)=V(m,0,0,\ldots)$ as in 1.2.2. Then equation (1.24) (with $a=1$) becomes
  $$ 0= \del_x\biggl(\frac F{\del_xF}\,\del_sF-\frac{xF^3}{(\del_xF)^2}\,\del_xF\biggr)-\del_s F 
    =\frac\del{\del x}\biggl[\frac xs\,f(xs)-\frac xs\,\frac{f(xs)^3}{f'(xs)}\biggr]-\frac\del{\del s}\,f(xs)$$
and if we write $f=y'$ then we see that this is (up to a power of $s$) the derivative of the differential equation $y=x{y'}^3/y''$ for the function $y=\sum\dfrac{v(n)x^{n-2}}{(n-2)!(n-3)!}$ discussed in the Introduction.

\newpage

\centerline{\bf \S \, 2. Explicit formulas and the inversion of the generating function.}
\bigskip

{\bf 2.1. Notation.} In this section we fix a value of the genus $g \geq 0$. We keep $g$ in the notation 
for $\mbgn$ and $V_g(\mb)$ but skip it elsewhere. To state our explicit formulas we must introduce some 
additional classes in $H^*(\mbgn, \Q)$. Recall that 
  $$ \Psi_{n,i} := \xi_i^* (c_1(\wcm)) \in H^2(\mbgn, \Q) \tag2.1$$
where $\xi_i: \mbgn \rightarrow \Cal C_n$ are the structure sections of the universal curve. 

After Witten [W], the integrals of top degree monomials in $\Psi_{n,i}$ are denoted
  $$ \langle\tau_{a_1}\dots\tau_{a_n}\rangle=\int_{\mbgn}\Psi_{n,1}^{a_1}\dots\Psi_{n,n}^{a_n}\tag 2.2 $$
Below we will express $V_g(\mb)$ via these numbers. For $g=0$,
they are just multinomial coefficients:
  $$\langle\tau_{a_1}\dots\tau_{a_n}\rangle_{g=0}=\frac{(a_1+\dots+a_n)!}{a_1!\dots a_n!} \tag2.3 $$
(see e.g. [K2], p.354). The structure of a generating series for all 
$\langle \tau_{a_1} \dots \tau_{a_n} \rangle$ and all $g$ was predicted by 
Witten [W], and Kontsevich identified it as a ``matrix Airy function'',  cf. [K1] and below.

More generally, we will consider the relative integrals of the type (2.2).
For $ k \geq l$, denote by $\pi_{k,l}: \Mb_{g,k} \rightarrow \Mb_{g,l}$ 
the morphism forgetting the last $k-l$ points. For any $a_1, \dots, a_p \geq 0$ define
  $$ \k_n(a_1, \dots, a_p) := \pi_{n+p,n*} (\Psi_{n+p,n+1}^{a_1+1} \dots \Psi_{n+p,n+p}^{a_p+1}) 
   \in H^{2(a_1+\dots +a_p)}(\mbgn,\Q). \tag 2.4 $$
Notice that whenever $a_1+\dots+a_p=\dim\mbgn$, we have also $(a_1+1)+\dots+(a_p+1)=\dim\Mb_{g,n+p}$, and then
  $$ \int_{\mbgn} \k_n(a_1,\dots, a_p) = \int_{\Mb_{g,n+p}} \Psi_{n+p,n+1}^{a_1+1} \dots \Psi_{n+p,n+p}^{a_p+1}
  = \langle \tau_0^n\tau_{a_1+1} \dots \tau_{a_p+1} \rangle .  \tag 2.5  $$

\proclaim {\quad 2.2. Theorem}  For any $g,n,a_1,\dots,a_p, a_i \ge 0$ we have
  $$ \w_n(a_1) \dots \w_n(a_p) = \sum_{k=1}^p \frac{(-1)^{p-k}}{k!} 
  \sum \Sb \{1,\dots,p\} = S_1 \amalg \dots \amalg S_k \\ S_i \neq \emptyset\endSb
  \k_n ( \sum_{j \in S_1} a_j, \dots, \sum_{j \in S_k} a_j) \tag 2.6 $$
Equivalently, for any $\mb \in N^\infty \setminus \{\ol{0}\}, p = \|\mb\|$,
 $$ \frac{(-1)^p}{\mb !}\bw_n^\mb  = \sum_{k=1}^p \frac {(-1)^{k}}{k!} 
 \sum \Sb \mb= \mb_1 + \dots + \mb_k \\ \mb_i \neq 0 \endSb
 \frac {\k_n(|\mb_1| , \dots , |\mb_k|)}{\mb_1!\dots \mb_k!} \tag 2.7 $$  \endproclaim
The proof consists of a geometric and a combinatorial part.

The first one is summarized in [AC], (1.12), (1.13) and given here in a slightly different notation. 
It was previously obtained by C\. Faber and D\. Zagier (unpublished).

\proclaim {\quad \bf 2.2.1. Lemma} We have  
  $$ \k_n(a_1, \dots, a_p) = \sum_{\s \in \Bbb S_p} 
  \prod_{\roman o\in \roman o(\s)} \w_n(\sum_{j\in \roman o} a_j) \tag 2.8 $$
where $\roman o(\s)$ denotes the set of the cycles of $\s$ acting on
$\{1,\dots,p\}$, i\.e\. the orbits of the cyclic group $\langle \s \rangle$.  \endproclaim

E. Albarello and M. Cornalba ([AC]) show that (2.8) formally follows with the help of the
push-pull formula from an identity going back to Witten [W]
  $$ \pi_{n+1,n*} (\Psi_{n+1,1}^{a_1} \dots \Psi_{n+1,n}^{a_n} \Psi_{n+1,n+1}^{a_{n+1}+1})
   = \Psi_{n,1}^{a_1} \dots \Psi_{n,n}^{a_n} \w_n(a_{n+1}) \tag 2.9 $$
and an identity for which a geometric proof is supplied in [AC]:
  $$ \w_n(a) = \pi_{n,n-1}^*(\w_{n-1}(a)) + \Psi_{n,n}^a. \tag 2.10 $$
We will not repeat their argument here.

The passage from (2.8) to (2.6) and (2.7) is a formal inversion result which we will prove here in an axiomatized
form because it is useful in other contexts as well.

Let $R$ be a commutative $\Q$-algebra generated by some elements $\w (a)$ where $a$ runs over all elements
of an additive semigroup $A.$

\proclaim{\quad 2.2.2. Lemma} Define elements $\w (a_1,\dots ,a_p)\in R$ for $p\ge 2,\ a_i\in A$ recursively by
 $$ \w(a_1,\dots,a_p)=\w(a_1,\dots,a_{p-1})\w(a_p)+\sum_{i=1}^{p-1}\w(a_1,\dots,a_i+a_p,\dots,a_{p-1}).\tag 2.11 $$
Then $\{\w (a_1,\dots ,a_p)\,|\,p\ge 1\}$ span $R$ as a linear space. They can be expressed via monomials
in $\w (a)$ by the following universal identity (coinciding with $(2.8)$):
 $$\w(a_1,\dots,a_p)=\sum_{\s\in\Bbb S_p}\prod_{\roman o\in\roman o(\s)}\w(\sum_{j\in\roman o}a_j)\tag2.12$$
In particular, $\w (a_1,\dots ,a_p)$ are symmetric in $a_1,\dots ,a_p.$

Conversely, monomials in $\w (a)$ can be expressed via
these elements by the universal formula (coinciding with $(2.6)$):
  $$ \w (a_1) \dots \w (a_p) = \sum_{k=1}^p \frac{(-1)^{p-k}}{k!} 
  \sum \Sb \{1,\dots,p\} = S_1 \amalg \dots \amalg S_k \\ S_i \neq \emptyset\endSb
  \w ( \sum_{j \in S_1} a_j, \dots, \sum_{j \in S_k} a_j). \tag 2.13 $$
If $A=\{ 1,2,3,\dots \}$, we have also
$$ \frac{(-1)^p}{\mb !}\bw_n^\mb  =  \sum_{k=1}^p \frac {(-1)^{k}}{k!} 
  \sum\Sb\mb=\mb_1+\dots+\mb_k\\ \mb_i\ne0\endSb \frac{\k_n(|\mb_1|,\dots,|\mb_k|)}{\mb_1!\dots\mb_k!}\tag 2.14 $$
where $p=\|\mb\|$ as in $(2.7)$ and $\bw^\mb=\w(1)^{m(1)}\w(2)^{m(2)}\dots$ .

Furthermore, $\w (a)$ are algebraically independent iff $\w (a_1,\dots ,a_p)$ are linearly independent.
In this case $R$ is graded by $A$ via $\roman{deg}\, \w (a_1,\dots ,a_p)=a_1+\dots +a_p.$ \endproclaim

\smallskip
{\bf Example.} The elements $\w (a_1,\dots ,a_p)$ are given for $p=2$ by
  $$ \w(a,b)=\w(a)\w(b)+\w(a+b),\quad \w(a)\w(b)=\w(a,b)-\w(a+b), $$
and for $p=3$ by
  $$ \align \w(a,b,c)&=\w(a)\w(b)\w(c)+\w(a+b)\w(c)+\w(a+c)\w(b)\\
               &\qquad\quad+\w(b+c)\w(a)+2\w(a+b+c), \\
      \w(a)\w(b)\w(c)&=\w(a,b,c)-\w(a+b,c)-w(a+c,b)\\
               &\qquad\quad-\w(b+c,a)+\w(a+b+c). \endalign $$
\smallskip
{\bf Proof of Lemma 2.2.2.} The following identity shows by induction in $p$ that (2.11) and (2.12) are equivalent:
$$
\multline \sum_{\s\in\Bbb S_{p+1}}\prod_{\roman o\in\roman o(\s)}\w(\sum_{j\in\roman o}a_j) = \\
\w(a_{p+1})\biggl[\sum_{\s\in\Bbb S_p}\prod_{\roman o\in\roman o(\s)}\w(\sum_{j\in\roman o}a_j)\biggr]
  +\sum_{i=1}^p\sum_{\tau\in\Bbb S_p}\prod_{\roman o\in \roman o(\tau)}
   \w(\sum_{j\in\roman o}a_j+\delta_{ij}a_{p+1})\endmultline \tag 2.15 $$
To convince yourself of the validity of (2.15) look at the following
bijective map from the \lhs monomials in $\w (a)$ to the \rhs monomials.
If a \lhs monomial is indexed by $\s \in \Bbb S_{p+1}$ for which 
$\s(p+1) = p+1$, we get it in the \rhs for $\tau =$ restriction of $\s$ to
$\{1, \dots,p\}$. Otherwise $p+1$ belongs to an orbit of $\s$ of cardinality
$\geq 2$; deleting $p+1$ from this cycle, we get a permutation $\tau \in \Bbb S_p$  and a number
$i = \s(p+1) \leq p$ producing exactly the needed monomial in the second sum.
 
One can similarly pass from (2.11) to (2.13) and backwards.
For example, assume that (2.13) is already proved for some $p$. Then we have
$$
\multline
[\w (a_1) \dots \w (a_p)]\w (a_{p+1}) = \\
\left ( \sum_{k=1}^p \frac{(-1)^{p-k}}{k!} 
\sum \Sb \{1,\dots,p\} = S_1 \amalg \dots \amalg S_k \\ S_i \neq \emptyset
\endSb
\k ( \sum_{j \in S_1} a_j, \dots, \sum_{j \in S_k} a_j) \right )
 \w (a_{p+1}) \\
= \sum_{k=1}^p \frac{(-1)^{p-k}}{k!} 
\sum \Sb \{1,\dots,p\} = S_1 \amalg \dots \amalg S_k \\ S_i \neq \emptyset
\endSb
\left ( \w ( \sum_{j \in S_1} a_j, \dots \sum_{j \in S_k} a_j,a_{p+1})- 
\right. \\
\sum_{i=1}^p \w ( \sum_{j \in S_1} (a_j+\delta_{ij}a_{p+1}, \dots ,
\left. \sum_{j \in S_k} (a_j+\delta_{ij}a_{p+1})) \right) .
\endmultline
\tag 2.16
$$
Now essentially the same combinatorics as above govern a correspondence between 
the summands in (2.16) and those in the \rhs of (2.13) for $p+1$ arguments which is
$$
\sum_{k=1}^{p+1} \frac{(-1)^{p+1-k}}{k!} 
\sum \Sb \{1,\dots,p+1\} = S_1 \amalg \dots \amalg S_k \\ S_i \neq \emptyset
\endSb
\w ( \sum_{j \in S_1} a_j, \dots, \sum_{j \in S_k} a_j).
$$
Namely, any ordered $k$-partition 
$\{1,\dots,p\} = S_1 \amalg \dots \amalg S_k $
can be enhanced to $k+1$ ordered $(k+1)$-partitions of
$\{1,\dots,p+1\}$ containing $\{p+1\}$ as a separate part, and to $k$ ordered
$k$-partitions of $\{1,\dots,p+1\}$ for which $p+1$ is put into one of the $S_i$'s.

It remains to rewrite (2.13) in the form (2.14), when
$A=\{ 1,2,3,\dots \}$. To this end, notice that if 
$\delta_{a_1}+ \dots + \delta_{a_p} = \mb$, we have
 $$ \w(a_1)\dots\w(a_p)=\prod_{a\ge1}\w(a)^{m(a)}=\bw^\mb,\qquad p=\|\mb\|$$
and $m(a) = \roman{card}\{j\,|\,a_j =  a\}$. Any set partition 
$\{1,\dots,p\}= S_1 \amalg \dots \amalg S_k, \, S_i \neq 0$, produces a vector partition
  $$ \mb=\mb_1+\dots+\mb_k,\quad \mb_i=(m_i(a))\ne0,\quad m_i(a)=\roman{card}\{j\,|\,a_j=a\}$$
We have $\sum_{j \in S_i} a_j = |\mb_i|$, and each $\mb = \mb_1 + \dots + \mb_k$ comes from 
  $$ \prod_{a \geq1} \frac {m(a)!}{m_1(a)! \dots m_k(a)!}=\frac {\mb!}{\mb_1! \dots \mb_k!} $$
set partitions. This finishes the proof of Lemma 2.2.2 and Theorem 2.2.
For a further discussion of this combinatorial setting, cf. 2.6 below.

As a corollary, we get:

\proclaim{\quad 2.3. Corollary} We have for $p= \|\mb\|,\ 3g-3+n=|\mb |$:
  $$ V_g(\mb) = \frac 1{|\mb|!} \sum_{k=1}^p \frac{(-1)^{p-k}}{k!}
  \sum \Sb \mb = \mb_1 + \dots + \mb_k \\ \mb_i \neq 0 \endSb
  \frac{ \langle \tau_0^n\tau_{|\mb_1|+1} \dots \tau_{|\mb_k|+1} \rangle} {\mb_1! \dots \mb_k!}\tag 2.17 $$
In particular, if $g=0$, then
  $$ V(\mb) = \sum_{k=1}^p (-1)^{p-k} \binom {|\mb|+k}{k}
  \sum \Sb \mb = \mb_1 + \dots + \mb_k \\ \mb_i \neq 0 \endSb
  \frac1{\prod_{k=1}^p(|\mb_i|+1)!\mb_i!}\tag 2.18 $$   \endproclaim

{\bf Proof.} Combine (2.6), (2.7), (2.5) and (2.3).

\smallskip

{\bf 2.3.1. A special case.}
Putting in (2.16) $\mb = (n-3,0,0,\dots)$ and multiplying by $(n-3)!^2$,
we get the following formula for Zograf's numbers (0.3):
$$
v_n = \sum_{k=1}^{n-3} (-1)^{n-3-k} \frac {(n-3+k)!(n-3)!}{k!}
\sum \Sb n-3 = m_1 + \dots + m_k \\ m_i \neq 0 \endSb
\frac 1{\prod_{i=1}^k (m_i+1)!\, m_i!}
$$
which is equivalent to (0.7).

We now proceed to the generalization of (0.8).

\proclaim{\quad 2.4. Theorem} In the ring of formal series of one variable with coefficients in 
$\Q[\bold s]=\Q[s_1,s_2, \dots]$ we have the following inversion formula
  $$ y = \sum_{|\mb| \geq 0} V(\mb)\,\frac {x^{|\mb|+1}}{|\mb|+1}\,\bold s^\mb \Longleftrightarrow
  x = \sum_{|\mb| \geq 0}\frac {y^{|\mb|+1}}{(|\mb|+1)!}\frac {(-\bold s)^\mb}{\mb!}. \tag 2.19 $$
\endproclaim 
There are two ways to prove this theorem, one starting from the explicit formula (2.18) and the other
 using the differential equation for the generating function $F(x;\bold s)$ derived in \S1. 
Since we do not know which proof, if either, may be generalizable to the higher genus case, we will give both.

\smallskip
{\bf 2.4.1. First proof: explicit formula.}  From (2.18) we have for any $\mu \geq 1$:
  $$ \align \sum_{\mb: |\mb| = \mu} V(\mb)\,\bold s^\mb
  & = \left. \sum_{k=1}^\infty (-1)^k \binom {\mu+k}{k} 
  \biggl(\sum_{|\mb| >0} \frac {(-\bold s)^\mb}{(|\mb|+1)!\mb!} \biggr)^k \right|_{\roman{degree}\,\mu} \\ 
  &=\left.\biggl(\sum_\mb\frac{(-\bold s)^\mb}{(|\mb|+1)!\mb!}\biggr)^{-\mu-1}\right|_{\roman{degree}\,\mu} \\
  &=\;\text{coeff\. of $y^{\mu}$ in } \biggl(\frac{x(y)}y\biggr)^{-\mu-1}\,,\endalign$$
where
 $$ x=x(y;\bold s):=\sum_{|\mb|\ge0}\frac{y^{|\mb|+1}}{(|\mb|+1)!}\frac{(-\bold s)^\mb}{\mb!}\in\Q[\bold s][[y]] $$
is the power series occurring on the right-hand side of (2.19) and have used the binomial identity 
$\sum_{k=1}^\infty (-1)^k \binom {\mu+k}{k} z^k = (1+z)^{-(\mu+1)}-1$.
But for any power series $x(y)= \sum_{r\geq 1} b_r y^r$, $b_1\ne0$, we have
  $$ \split \text{coeff\. of $y^{\mu}$ in } \fracwithdelims() {x(y)}{y}^{-\mu-1}
  &= \roman {res}_{y=0} \biggl( \frac 1{y^{\mu+1}} \fracwithdelims() {y}{x(y)}^{\mu+1} dy\biggr) \\
  &= \roman {res}_{x=0} \biggl(\frac 1{x^{\mu+1}}\,\frac{dy(x)}{dx}\, dx\biggr) \\
  &= \text {coeff\. of $x^{\mu}$ in } \frac{dy(x)}{dx} \\
  &= \text {coeff\. of $\frac{x^{\mu+1}}{\mu+1}$ in } y(x)\;, \endsplit $$
where $y(x)$ is the power series obtained by formal inversion of $x = x(y)$ (which is
possible since $b_1 \neq 0$). Applying this to our case, we find that the inverse series of $x(y)$ is given by
 $$ y(x) = \sum_{\mu\ge0}\frac{x^{\mu+1}}{\mu+1}\sum_{|\mb|=\mu}V(\mb)\,\bold s^\mb\,,$$
which is precisely the expression on the left of (2.19).
\smallskip

{\bf 2.4.2. Second proof: differential equation.} 
In 1.6 we characterized the generating series (1.23) as the unique power series $F(x;\bold s)\in1+x\Q[\bold s][[x]]$
with $F(x;\bold0)=0$ for which there are power series $H_a(x,\bold s)$ satisfying
  $$ H_0(x)=x,\;\; H_a(0)=0,\;\;\del_x(H_a)=H_{a-1}F,\;\;\del_aF=\del_x(H_aF)\qquad(a\ge1)\,.\tag 2.20$$
Write $y=\int_0^x F(\xi;\bold s) d \xi=x+\cdots$ for the integral of $F$ with respect to $x$. 
This is an invertible power series, so
we can also write $x=x(y)=y+\cdots$ and define power series $h_a(y)=h_a(y;\bold s)\in y\Q[\bold s][[y]]$
by $h_a(y)=H_a(x(y))$. In terms of $h_a$, the first three equations in (2.20) become
  $$ h_0(y)=x(y),\quad h_0(0)=0,\quad h_a'(y)=h_{a-1}(y)\tag 2.21$$
while the last equation in (2.20) becomes
  $$ \del_a x(y,\bold s) = -h_a(y)\qquad(a\ge1).\tag 2.22$$
(To see this, first integrate the last equation in (2.20) to get $\del_ay=FH_a$, where the meaning of 
$\del_ay$ is that we differentiate in $s_a$ keeping $x$ constant.  But differentiating the identity
$x(y(x;\bold s);\bold s)\equiv x$ with resepct to $s_a$ gives $\del_ax+\del_ay\cdot\del_yx=0$, where
$\del_ax$ is the derivative of $x$ with respect to $s_a$ keeping $y$ constant, so
$\del_a x = -F^{-1}\del_ay=-H_a(x)=-h_a(y)$ as claimed.)

Equations (2.21) and (2.22) form a system of linear differential equations replacing the non-linear system (2.20). 
They can be combined into a system of linear differential equations for the single function $x(y)$, namely
 $$\frac{\del^2x}{\del s_1\del y}=-x\,,\qquad
  \frac{\del^2x}{\del s_a\del y}=\frac{\del x}{\del s_{a-1}}\qquad(a\ge2).\tag 2.23$$
But in fact this is not needed because we can solve the system immediately.  Write
 $$x=x(y)=x(y,\bold s)=\sum_{n=0}^\infty B_n(\bold s)\,\frac{y^{n+1}}{(n+1)!}\tag 2.24$$
with $B_0(\bold s)=1$ and $B_n(\bold 0)=0$ for all $n\ge1$. The solution of (2.21) is immediately seen to be
 $$h_a(y,\bold s)=\sum_{n=0}^\infty B_n(\bold s)\,\frac{y^{n+a+1}}{(n+a+1)!}\,,$$
and (2.22) then says that $\del_aB_n=-B_{n-a}$ ($=0$ if $a>n$), giving successively
$B_a=-s_1$, $B_2=\frac12s_a^2-s_2$, \dots and in general
 $$B_n(s)=\sum_{m_1+2m_2+\cdots=n}
  (-1)^{m_1+m_2+\cdots}\,\frac{\,s_1^{m_1}\,s_2^{m_2}\,\cdots\,}{m_1!\,m_2!\,\cdots}\;.\tag2.25$$
Substituting (2.25) into (2.24) gives the expansion on the \rhs of (2.19).

{\bf 2.4.3. Remarks.} Equation (2.19) specializes to (0.8) if we put $\bold s=(1,0,0,\dots)$, in which case
 the \rhs of (2.19) is a Bessel function.  In general, this \rhs cannot be expressed in terms of standard
functions, but the series is easily summed after applying to it the formal Laplace transform, since from
(2.24) and (2.25) we get immediately
  $$\eta^{-2}\,\int_0^\infty e^{-y/\eta}\,x(y)\,dy = e^{-s_1\eta-s_2\eta^2-\dots} \,. \tag 2.26 $$
Conversely, by integrating once by parts and making the change of variables from $y$ to $x=x(y)$, we find the dual formula
$$ \eta^{-1}\,\int_0^\infty e^{-y(x)/\eta}\,dx =  e^{-s_1\eta-s_2\eta^2-\dots}\,.\tag2.27$$
These formulas allow us to get analytic information about the generating series $y(x)$ (and hence about the
numbers $V(\mb)$ when we specialize $(s_1,s_2, \dots)$ suitably. 
A theoretical reinterpretation of them will be given in Theorem 3.4.2 below.

It would be interesting to extend (2.17) including the contributions of all genera and eventually of all 
combinatorial cohomology classes as in [K1]. Below we collect some remarks in this direction.

\smallskip

{\bf 2.5. Kontsevich's formulas.} We start with Kontsevich's formula ([K1], sec. 3.1) for the correlation numbers
$\langle \tau_{d_1}\dots\tau_{d_n}\rangle$ which can be used to calculate algorithmically $V_g(\mb )$ 
with the help of (2.17).  It has the following structure.  \def\l{\lambda}

Fix $(g,n),$ put $d=3g-3+n$, and choose $n$ independent variables $\l_1,\dots ,\l_n.$ Then
  $$ \sum \Sb d= d_1 + \dots + d_n \\ d_i \ge 0 \endSb \langle \tau_{d_1}\dots\tau_{d_n}\rangle
     \prod_{i=1}^n \frac{(2d_i-1)!!}{\l_i^{2d_i+1}}=
    \sum_{\Gamma\in G_{g,n}}\frac{2^{-|V_{\tau}|}}{|\roman{Aut}\,\tau|}
    \prod_{e\in E_{\tau}}\frac{2}{\l^\prime (e)+\l'' (e)}. \tag 2.28 $$ 
Here $G_{g,n}$ is the set of the isomorphism classes of triples $\Gamma =(\tau ,c,f)$ where: \smallskip

i) $\tau$ is a connected graph with all vertices $v\in V_{\tau}$ of valency 3 and no tails;

ii) $c$ is a family of cyclic orders on all $F_{\tau}(v)$ where $F_{\tau}(v)$ is the set of flags adjoining $v$;

iii) $f$ is a bijection between $\{ 1,\dots ,n\}$ and the set of all cycles of $\tau$. We recall that a cycle is a
cyclically ordered sequence of edges (without repetitions) $(e_1,e_2,\dots ,e_k)$ such that for every $i$, 
$e_i$ and $e_{i+1}$ have a common vertex $v_i$ and the flag $(e_i,v_i)$ follows the
flag $(e_{i+1},v_i)$ in the sense of $c$;

iv) for any edge $e\in E_{\tau}$, $\{\l^\prime (e),\l''(e)\}=\{\l_a,\l_b\},$ where 
$\{a,b\}\subset\{1,\dots ,n\}$ are the $f$-labels of the two cycles to which $e$ belongs.  \smallskip

If $\tau$ is embedded into a closed Riemann surface $X$ oriented compatibly with $c$, the cycles of $\tau$
become the boundaries of the oriented connected components of $X\setminus |\tau |$ (2-cells).
Then $f$ labels these cells, and $\{ a,b\}$ become the labels of the cells adjoining $e.$  \smallskip

A paradoxical property of (2.28) which does not allow to read off $\langle \tau_{d_1}\dots\tau_{d_n}\rangle$
directly from this identity is the cancellation of poles at $\l_a=-\l_b$ in the r\.h\.s\., not at all
evident a priori even in the simplest case $g=0,\,n=3:$
$$\align\frac{\langle\tau_0\tau_0\tau_0\rangle}{\l_1\l_2\l_3}&=\frac1{\l_1\l_2\l_3}
  =\frac1{\l_1(\l_1+\l_2)(\l_1+\l_3)}+\frac1{\l_2(\l_2+\l_3)(\l_2+\l_1)}\\
  &\qquad+\frac1{\l_3(\l_3+\l_1)(\l_3+\l_2)}+\frac1{(\l_1+\l_2)(\l_2+\l_3)(\l_3+\l_1)}. \endalign$$

A generating function incorporating all stable $(g,n)$
that can be summed using (2.28) and the standard techniques of perturbation theory is
 $$  W(s_0,s_1,\dots )=
  \langle \roman{exp}\,(\sum_{a=0}^\infty s_a\tau_a)\rangle=
  \sum_{n_i\ge 0}\langle \tau_0^{n_0}\tau_1^{n_1}\dots\rangle
  \frac{s_0^{n_0}s_1^{n_1}\dots}{n_0!\,n_1!\dots} \tag 2.29 $$

Kontsevich's theorem states that (2.29) is an asymptotic expansion of (the logarithm of) a matrix Airy function.

We will skip the description of this function because we were unable to find a sensible operator
processing $W(s_0,s_1,\dots )$ into a generating series for the WP-volumes.

Instead we will show that the formalism of Lemma 2.2.2 has a nice self-reproducing property in the language
of formal series, but in order to use it in our context,
a different generating series for $\langle \tau_{d_1}\dots\tau_{d_n}\rangle$ is needed.
\smallskip

{\bf 2.6. A remark on Lemma 2.2.2.} In the situation
of this Lemma, assume that $A=\{0,1,2,\dots \}$ and put for $p\ge 1$
  $$ U(t_1,\dots ,t_p):=\sum_{a_1,\dots ,a_p\ge 0}
  \w(a_1,\dots,a_p)\frac{t_1^{a_1}}{a_1!}\dots\frac{t_p^{a_p}}{a_p!}. \tag 2.30 $$

\proclaim{\quad 2.6.1. Lemma } We have
  $$ U(t_1,\dots ,t_p)=U(t_1,\dots ,t_{p-1})U(t_p)+\sum_{i=1}^{p-1} U(t_1,\dots ,t_i+t_p,\dots ,t_{p-1}). \tag 2.31 $$
\endproclaim
The proof is a straightforward calculation using (2.11).

\proclaim{\quad 2.6.2. Corollary} We have
  $$ U(t_1)\dots U(t_p)=\sum_{k=1}^{p}\frac{(-1)^{p-k}}{k!}
  \sum \Sb \{1,\dots,p\} = S_1 \amalg \dots \amalg S_k \\ S_i \neq \emptyset \endSb
  U(\sum_{j\in S_1}t_j,\dots ,\sum_{j\in S_p}t_j). \tag 2.32 $$
\endproclaim
\smallskip

For the proof, apply Lemma 2.22 to $A=\oplus_{i=1}^p\Bbb{Z}t_i,\, R[[t_i]],$ and
$U(t)$ instead of $\w (a).$  \smallskip

We can try to use (2.32) in the case $R=H^*(\mbgn),\w (a)=\w_{g,n}(a).$ The l\.h\.s\. of (2.32)
after integration becomes a polynomial with coefficients which are WP-volumes multiplied by some factorials,
whereas the r\.h\.s\. becomes a similar polynomial which is a linear combination of correlation numbers.

\newpage
\centerline{\bf \S 3. Potential of the invertible Cohomological Field Theories}
\bigskip

Here we explain following [KMK] that the genus zero generating function
for higher WP-volumes is the third derivative of the potential of a generic
invertible CohFT written in coordinates additive with respect to the
tensor multiplication.  \smallskip

\proclaim{\quad 3.1. Definition}
Let $k$ be a ground field of characteristic zero, $H$ a finite dimensional $\Bbb Z_2$-graded linear space over $k$, 
and $g$ an even nondegenerate symmetric pairing on $H$ represented by its dual element $\Delta \in H \otimes H$.

A structure of the tree level Cohomological Field Theory on $(H,g)$ is given by a sequence of maps
  $$ I_n: H^{\otimes n} \rightarrow H^* (\mbn,k), \qquad n \geq 3 $$
which are $\Bbb S_n$-covariant and satisfy the following identities.  
For any 2-partition $\s: \{1,\dots,n\} = S_1 \amalg S_2$,  $n_i = |S_i| \geq 2$ and the respective embedding
of the boundary divisor $\varphi_\s: \Mb_{0,n_1+1} \times \Mb_{0,n_2+1} \rightarrow \Mb_{0,n}$ we have
  $$ \varphi^*_\s (I_n(\gamma_1 \otimes \dots \otimes \gamma_n)) = 
  \epsilon(\s)(I_{n_1+1}\otimes I_{n_2+1})
  (\bigotimes_{j\in S_1}\gamma_j\otimes\Delta\otimes(\bigotimes_{j\in S_2}\gamma_j)) $$
where $\epsilon(\s)$ is the sign of $\s$ acting on the odd $\gamma_i \in H$.  \endproclaim

This notion was introduced in [KM]. One of its most interesting instances is 
quantum cohomology: a canonical structure of CohFT on  $(H^*(V,k),$ {\it Poincar\'e pairing})
for any smooth projective (or $C^\infty$ compact symplectic) manifold $V$.
\smallskip

{\bf 3.2. Tensor product.} Let $\{H', g', I'_n\}$ and $\{H'',g'',I''_n\}$ be two CohFT's. 
Put $H = H' \otimes H''$ and $g =g' \otimes g''$. We can define a CohFT on $(H,g)$ by
  $$ I_n(\gamma'_1 \otimes \gamma''_1 \otimes \dots \otimes \gamma'_n \otimes \gamma''_n) :=\epsilon(\gamma',\,\gamma'')
   I'_n(\gamma'_1\otimes\dots\otimes\gamma'_n)\wedge I''_n(\gamma''_1\otimes\dots\otimes\gamma''_n). $$
In the context of quantum cohomology, this product serves as a K\"unneth formula.
\smallskip

{\bf 3.3. Potential.}
 The potential of a CohFT $(H,g,I_n)$ is the formal series $\Phi \in k[[H^{\vee}]]$ which in terms of coordinates 
w\.r\.t\. a basis $(\Delta_a)$ of $H$ can be written as
  $$ \Phi(x) := \sum_{n\geq 3}\frac 1{n!}\int_{\mbn}I_n((\sum x^a\Delta_a)^{\otimes n}). \tag 3.1 $$
The main result about $\Phi(x)$ proved in [KM] and [KMK] is the following theorem.

\proclaim{\quad 3.3.1. Theorem} 
The map $(H,g,I_n) \mapsto \Phi$ establishes a bijection between the following objects:  \roster
\item"a)" The structures of a CohFT on $(H,g)$
\item"b)" The solutions of the Associativity Equations in $k[[H^{\vee}]]$ modulo terms of degree $\le 2.$
\endroster \endproclaim
\smallskip

We recall that the Associativity Equations are
  $$ \forall a,b,c,d:\sum_{ef}\del_a\del_b\del_e\Phi\cdot g^{ef}\del_f\del_c\del_d\Phi
  =(-1)^{\tilde x_a(\tilde x_b+\tilde x_c)}\sum_{ef}\del_b\del_c\del_e\Phi\cdot g^{ef}\del_f\del_a\del_d\Phi $$
where $(g^{ef})=(g(\Delta_e,\Delta_f))^{-1},\del_a=\del/\del x^a,\tilde x^a=$ 
the $\Bbb Z_2$-parity of $x^a$ and $\Delta_a$.

\smallskip

{\bf 3.4. The moduli space of one-dimensional CohFT's.}
We can think naively of CohFT's as forming an infinite dimensional algebraic
variety, with the tensor product defining a structure of a semigroup on it. 
In view of 3.3.1, it is natural to try to understand the properties of the
potential as a function on the moduli space.  In particular, we would like to
understand how to express the potential function $\Phi_{\A'\otimes\A''}$ associated to the 
tensor product of two CohFT's $\A'=(H',g',I')$ and $\A''=(H'',g'',I'')$ in terms of the
potential functions  $\Phi_{\A'}$ and $\Phi_{\A''}$.

As a special case, let us consider CohFT structures on one-dimensional spaces. Such a theory will be invertible 
with respect to the tensor product if the map $I_3$ from $H^{\otimes3}\cong k$ to $H^*(\ol M_{03},k)\cong k$
is an isomorphism. We will call the theory normalized if we have a chosen
basis of length one, $H= k \Delta_0$, $g(\Delta_0,\Delta_0)=1$, and $I_3(\Delta_0\otimes\Delta_0\otimes\Delta_0)=1$.
(or equivalently $I_n(\Delta_0^{\otimes n})=1_n\,+$ terms of dimension $>0$
for all $n\geq3$, where $1_n\in H^0(\mbn,k)$ is the fundamental class.)  The potential function (3.1) has the form
$$ \Phi_\A(x)=\sum_{n=3}^\infty C_n\,\frac{x^n}{n!} \tag 3.2$$
where $C_3=1$ but the other coefficients are arbitrary by virtue of Theorem 3.3.1, since the Associativity Equations
are empty in this case.  Thus the space $\NIC$ of all normalized and invertible 1-dimensional CohFT's is
canonically isomorphic to $\frac16x^3\,+\,x^4\,k[[x]]$ and has canonical coordinates $C_n$ ($n\ge4$), and 
we would like to describe the tensor product in terms of these coordinates.

These 1-dimensional CohFT structures were studied in [KMK] and a different set of coordinates was given.
For $s_1,\,s_2,\dots\in k$ there is an element $\A(\bold s)\in\NIC$ given by
  $$I_n(\Delta_0^{\otimes n})=\w_n[s_1,s_2,\dots]:=\exp\bigl(\sum_{a=1}^\infty s_a\w_n(a)\bigr)\qquad(n\ge3),\tag 3.3 $$
in the notation of our paper and [KMK]. Then it was shown that the map $\bold s\mapsto \A(\bold s)$ gives a bijection
between $k^\N$ and $\NIC$ and that $\A(\bold s')\otimes\A(\bold s'')\cong\A(\bold s'+\bold s'')$, i.e.:
\proclaim{\quad 3.4.1. Theorem [KMK]}  The parameters $(s_1,s_2,\dots)$ form a coordinate system on the space of
normalized 1-dimensional CohFT's. The tensor product becomes addition in these coordinates. \endproclaim

Denote by $\Phi(x;\bold s)$ the potential associated to the theory (3.3).  The connection with what we have 
done in this paper is that the third derivative of the potential $\Phi(x;\bold s)$ associated to the 
theory (3.3) is just our generating function for higher WP-volumes.  Indeed, definition (3.1) gives
  $$ \Phi(x;s_1,s_2,\dots) = \sum_{n=3}^\infty \frac{x^n}{n!}\int_{\mbn}
   \sum_{\sum am(a)=n-3}\prod_{a\ge1}\w_n(a)^{m(a)}\frac{s_a^{m(a)}}{m(a)!}\,, $$
and the third derivative of this is obviously the function $F(x;\bold s)=\sum_\mb V(\mb)x^{|\mb|}\bold s^\mb$
defined in 1.6.  We now use this connection to describe both the tensor product and the coordinates on
the space of invertible 1-dimensional CohFT's explicitly.

\proclaim{\quad 3.4.2. Theorem}  Define bijections
$$ \NIC\LR\frac{x^3}6+x^4\,k[[x]]\LR 1+\eta\,k[[\eta]] \,,\tag 3.4$$
where the first map assigns to a theory $\A$ its potential $\Phi_\A(x)$ and the second map is defined by
$$\Phi(x)\,\leftrightarrow\, U(\eta)=\int_0^\infty e^{-\Phi''(\eta x)/\eta}\,dx\tag3.5$$
or alternatively by assigning to $\Phi(x)=\frac16x^3+\dots$ the power series $U(\eta)=\sum_{n=0}^\infty B_n\eta^n$
where $x=\sum B_n\frac{y^{n+1}}{(n+1)!}=y+\cdots$ is the inverse power series of $y=\Phi''(x)=x+\cdots$. 
Then the tensor product of 1-dimensional CohFT's corresponds to multiplication in $1+\eta k[[\eta]]\,$:
$U_{\A'\otimes\A''}(\eta)=U_{\A'}(\eta)\,U_{\A''}(\eta)$. The coefficients of $-\log U_\A(\eta)$ are
the canonical coordinates of $\A\,$.\endproclaim
{\bf Proof.}  Equations (2.25) and (2.26) tell us that the two descriptions of $U(\eta)$ given in the theorem 
agree and give $e^{-s_1\eta-s_2\eta^2-\cdots}$ when $\Phi(x)=\Phi(x;\bold s)$, and 
in view of Theorem 3.4.1 this implies the result in general. 
\smallskip
{\bf 3.5. Explicit formulas.} 
Substituting (3.2) (with $C_3=1$) into (3.5), expanding, and integrating term by term gives the explicit formula
 $$ B_n=\sum\Sb n_4,n_5,\ldots\ge0\\ n_4+2n_5+\cdots=n\endSb
 \frac{(2n_4+3n_5+\cdots)!}{2!^{n_4}\,3!^{n_5}\cdots\,n_4!\,n_5!\,\cdots}\,(-C_4)^{n_4}\,(-C_5)^{n_5}\cdots$$
for the coefficients of $U(\eta)$ in terms of the coefficients of $\Phi(x)$, and the same argument applied 
to the inverse power series gives the reciprocal formula
 $$ C_n=\sum\Sb n_1,n_2,\ldots\ge0\\ n_4+2n_5+\cdots=n-3\endSb
  \frac{(2n_1+3n_2+\cdots)!}{2!^{n_1}\,3!^{n_2}\cdots\,n_1!\,n_2!\,\cdots}\,(-B_1)^{n_1}\,(-B_2)^{n_2}\cdots\,.$$
Combining these formulas with the identity $U_{\A'\otimes\A''}(\eta)=U_{\A'}(\eta)\,U_{\A''}(\eta)$ we obtain
the explicit law for the tensor product of two normalized invertible CohFT's in terms of the coefficients of their
potential functions:
$$ \align C_4 &=C_4'+C_4''\,,\\  C_5 &=C_5'+5C_4'C_4''+C_5''\,,\\ 
   C_6 &=C_6'+(8\,{C_4'}^2+C_5')\,C_4''+C_4'\,(8\,{C_4''}^2+C_5'')+C_6''\,,\\
   C_7 &=C_7'+ (35\,C_4'\,C_5'+14\,C_6')\,C_4''+(61\,{C_4'}^2\,{C_4''}^2+33\,{C_4'}^2C_5''+33\,C_5'\,{C_4''}^2 \\
   &\qquad+19\,C_5'\,C_5'') + C_4'\,(35\,C_4''\,C_5''+14\,C_6'') +C_7''\,,\qquad\ldots \endalign$$

Finally, we observe that the values of the genus 0 Weil-Petersson volumes $V(\mb)$ can be calculated numerically
from any of a number of formulas in this paper: the recursion relation (1.4), the differential equation (1.24),
the closed formula (2.18), or the generating function formula (2.19). Here are the
values up to $|\mb|=5$, expressed in terms of the generating function (1.23):
  $$ \align F(&x,\bold s) = 1\,+\,s_1\,x\,+\,\bigl(5\,\frac{s_1^2}2+s_1\bigr)\,\frac{x^2}2
\,+\,\bigl(61\,\frac{s_1^3}6+9\,s_1s_2+s_3\bigr)\,\frac{x^3}6\\
&+\, \bigl(1379\,\frac{s_1^4}{24}+161\,\frac{s_1^2s_2}2+14\,s_1s_3+19\,\frac{s_2^2}2+s_4\bigr)\,\frac{x^4}{24} \\
&+\, \bigl(49946\,\frac{s_1^5}{120}+4822\,\frac{s_1^3s_2}6+344\,\frac{s_1^2s_3}2+470\,\frac{s_1s_3^2}2+20\,s_1s_4
    +34\,s_2s_3+s_5\bigr)\,\frac{x^5}{120} \\
&+\,\text O(x^6)\,. \endalign$$
Note that the coefficient $\dsize\int_{\ol M_{0,n}}\bw^\mb$ of $\dfrac{\bold s^\mb}{\mb!}\,\dfrac{x^{|\mb|}}{|\mb|!}$
is integral for every $\mb$ since the cohomology classes $\wgn(a)$ are integral for $g=0$.

\newpage

\centerline{\bf Appendix}

\bigskip
In this appendix we make a few remarks on the asymptotics of the Betti numbers and Euler characteristics of 
the moduli spaces of $\mbn$ and on the Weil-Petersson volumes(0.2). Set $P_1(q)=1$ and for $n\ge2$ let 
  $$B_j(n)=\dim H^j(\overline{M}_{0,n+1})\,,\qquad P_n(q)=\sum_{j\ge0}B_j(n)\,q^j$$
be the Betti numbers and Poincar\'e polynomial, respectively, of $\overline{M}_{0,n+1}$.
It was shown in [M] that the polynomials $P_n$ satisfy the recursion
  $$ P_{n+1}(q)=P_n(q)\,+\,q^2\sum_{m=2}^n\binom nmP_m(q)\,P_{n+1-m}(q)\qquad(n\ge1)\,.\tag A.1$$
This is equivalent to the differential equation 
  $$ \frac{\del y}{\del x}=\frac{1+y}{1-q^2(y-x)}\,,\tag A.2$$
for the generating function
  $$y =\sum_{n=1}^\infty P_n(q)\,\frac{x^n}{n!}=x+\frac12\,x^2+\frac{1+q^2}6\,x^3 +\frac{1+5q^2+q^4}{24}\,x^4+\cdots$$
The solution is given by the implicit equation
  $$(1+y)^{q^2}\,=\,1\,+\,q^2\,x\,+\,q^4\,(y-x)\,.\tag A.3$$

If we specialize to $q=1$, then the solution of (A.2), or the limiting value of (A.3), is given by the implicit equation 
  $$x=2y-(1+y)\log(1+y)=y-\sum_{n=2}^\infty\frac{(-1)^ny^n}{n(n-1)}\qquad(q=1)\,,$$
from which the asymptotic formula (0.10) for the values $\chi(\overline{M}_{0,n+1})=P_n(1)$ follows
easily. (The derivative $dx/dy$ vanishes simply  at $y=e-1$, $x=e-2$, so the power series expansion
of $y(1,x)$ as a function of $x$ has radius of convergence $e-2$ with a square-root singularity
at $x=e-2$.) The same method applies to the inversion formula (0.8) to give the asymptotic equation
(0.9) with the constant $C$ given by $C=2\gamma_0J_0'(\gamma_0)$ where $\gamma_0=2.4048255577$ is the 
smallest zero of the Bessel function $J_0(x)$ and the other constants in the expansion can be obtained 
by doing a more detailed analysis of the function $\sqrt y\,J_1(2\sqrt y)$ near its maximum.

We can also use (A.1) and (A.3) to study the behavior of the Betti numbers $B_j(n)$ as 
a function of $n$ for fixed~$j$.  From (A.1) we get
  $$\align B_2(n) &=2^n-\frac12(n^2+n+2),\\ 
  B_4(n) &=\frac32\,3^n-\frac14(n^2+5n+8)\,2^n+\frac1{24}(3n^4+2n^3+21n^2+22n+12)\endalign$$
and more generally 
  $$ B_{2j}(n)=\sum_{l=0}^j p_{j,l}(n)\,(j+1-l)^n\tag A.4$$
for some polynomials $p_{j,l}(n)$ of degree $2l$ in $n$. (The odd Betti numbers are 0.) To
see this, and to get information about the polynomials $p_{j,l}$, we observe that (A.4) is
equivalent to the statement that for each $j$ we have the generating function identity
  $$ \sum_{n=1}^\infty B_j(n)\,\frac{x^n}{n!}=A_j(x,\,e^x)$$
for some polynomial $A_j(x,u)\in\Bbb Q[u,x]$, the first few values being
  $$\align & A_0=u-1,\qquad A_2=u^2-\bigl(1+x+\frac{x^2}2\bigr)\,u,\\ 
 & A_4=\frac32\,u^3-(2+3x+x^2)\,u^2+\bigl(\frac12+2x+2x^2+\frac{5x^3}6+\frac{x^4}8\bigr)\,u\,.\endalign$$
The function $y$ is then replaced by a power series $y(q,x,u)=\sum A_j(x,u)\,q^j$ in three variables,
the previous power series being obtained by setting $u=e^x$, and the recursion (A.1) gives the same
differential equation (A.2) as before but with the left-hand side replaced by 
$\dfrac{\del y}{\del x}+u\,\dfrac{\del y}{\del u}$.  The solution is given by the implicit equation
  $$\biggl(\frac{1+y}{ue^{-x}}\biggr)^{q^2}\,=\,1\,+\,q^2\,x\,+\,q^4\,(y-x)\,.\tag A.5$$
The polynomial $A_j(u,x)$ has degree $j+2$, where $u$ and $x$ are assigned degrees 2 and 1, respectively, and
$u$ divides $A_j$ for $j\ge1$, so $q^2(1+y)=\Phi(q,qx,q^2u)$ with $\Phi(q,X,U)\in\Bbb Q[[q,X,U]]$. Then (A.5) gives  $$\Phi=U\,e^{-X/q}\,\bigl(1+qX+q^2\,\Phi+\text O(q^3)\bigr)^{1/q^2}=U\,e^{\Phi-X^2/2+\text O(q)}\,,$$
so $\Phi(0,X,U)$ is a solution of $\Phi e^{-\Phi}=Ue^{-X^2/2}$.  Inverting this gives
 $$1+y = q^{-2}\,\Phi=q^{-2}\,\sum_{j=1}^\infty\frac{(j+1)^{j-1}}{j!}\,U^{j+1}\,e^{-(j+1)X^2/2}\,(1+\text O(q)),$$
which with $u=Uq^{-2}$, $x=Xq^{-1}$ translates back as
  $$ A_j(x,u)=\sum_{l=0}^j\frac{(-1)^l}{2^ll!}\,\frac{(j-l+1)^{j-1}}{(j-l)!}\,u^{j-l+1}\,(x^{2l}+\text O(x^{2l-1}))$$
or, for the polynomials $p_{j,l}$ defined in (A.4), as
  $$ p_{j,l}(n) = \frac{(-1)^l}{2^ll!}\,\frac{(j-l+1)^{j-2l-1}}{(j-l)!}\,n^{2l}\,+\,\text O(n^{2l-1})\,.$$
In particular, 
  $$ B_{2j}(n)=\frac{(j+1)^{n+j-1}}{j!} + \text O\bigl(n^2\,j^n\bigr)\,.$$

\newpage

\centerline{\bf References}

\bigskip

[AC] E\. Arbarello, M\. Cornalba. {\it Combinatorial and algebro-geometric
cohomology classes on the moduli spaces of curves.} Preprint, 1995 
(to appear in J. of Alg. Geom.)

\smallskip

[G] E\. Getzler. {\it Operads and moduli spaces of genus zero
Riemann surfaces.} In: The Moduli Space of Curves, ed. by
R\. Dijkgraaf, C\. Faber, G\. van der Geer, Progress in Math\.
vol\. 129, Birkh\"auser, 1995, 199--230.

\smallskip

[Ke] S\. Keel. {\it Intersection theory of moduli spaces of stable
$n$-pointed curves of genus zero.} Trans. AMS, 330 (1992), 545--574.

\smallskip

[K1] M\. Kontsevich. {\it Intersection theory on the moduli space of curves and
the matrix Airy function.} Comm\. Math\. Phys\. 147 (1992), 1--23.

\smallskip

[K2] M\. Kontsevich. {\it Enumeration of rational curves via torus actions.}
In: The Moduli Space of Curves, ed. by
R\. Dijkgraaf, C\. Faber, G\. van der Geer, Progress in Math\.
vol\. 129, Birkh\"auser, 1995, 335--368.

\smallskip

[KM] M\. Kontsevich, Yu\. Manin. {\it Gromov-Witten classes, quantum
cohomology, and enumerative geometry.} Comm. Math. Phys.,
164:3 (1994), 525--562.

\smallskip

[KMK] M\. Kontsevich, Yu\. Manin (with Appendix by R\. Kaufmann).
{\it Quantum cohomology of a product.} Inv. Math., 124 (1996),
f. 1--3, 313--340.

\smallskip

[LS] B\. Logan, L\. Shepp. 
{\it A variational problem for random Young tableaux.} Adv\. in Math\., 
26 (1977), 206--222.

\smallskip

[M] Yu\. Manin. {\it Generating functions in algebraic geometry and sums over
trees.} In: The Moduli Space of Curves, ed. by
R\. Dijkgraaf, C\. Faber, G\. van der Geer, Progress in Math\.
vol\. 129, Birkh\"auser, 1995, 401--418.

\smallskip

[Ma] M\. Matone. {\it Nonperturbative model of Liouville gravity.}  Preprint hep-th/9402081.

\smallskip

[W] E\. Witten. {\it Two-dimensional gravity and intersection theory on 
moduli space.} Surv\. in Diff\. Geo\. 1 (1991), 243--310.

\smallskip

[VK] A. Vershik, S. Kerov. {\it Asymptotic theory of the characters of the
symmetric group.} Func. An. Appl., 15:4 (1981), 15--27.

\smallskip

[Z] P. Zograf. {\it The Weil-Petersson volume of the moduli spaces
of punctured spheres.} In: Cont\. Math\., 150 (1993), ed. by R\. M\. Hain and 
C\. F\. B\"odigheimer, 267--372.

\enddocument